\documentclass[preprint,12pt]{elsarticle}

\usepackage{graphicx}

\usepackage{amssymb}
\usepackage{amsmath}
\usepackage{amsthm}
\usepackage{subfigure}
\usepackage{subfiles}
\usepackage{arydshln}

\newtheorem{theorem}{Theorem}[section]
\newtheorem{lemma}{Lemma}[section]
\newtheorem{proposition}{Proposition}[section]
\newtheorem{assumption}{Assumption}[section]

\theoremstyle{definition}

\newtheorem{remark}{Remark}[section]
\def\dfrac{\displaystyle\frac}

\journal{Journal of the Franklin Institute}

\begin{document}

\begin{frontmatter}

\title{Output regulation for a reaction-diffusion system with input delay and unknown frequency}

\author[s1]{Shen Wang}
\ead{wangss@tju.edu.cn}

\author[s1]{Zhong-Jie Han\corref{1}}
\cortext[1]{Corresponding author}
\ead{zjhan@tju.edu.cn} 

\author[s1]{Kai Liu}
\ead{k\_liu@tju.edu.cn} 

\affiliation[s1]{
organization={School of Mathematics and KL-AAGDM, Tianjin University},
city={Tianjin},
postcode={300350}, 
country={China}}

\author[s2,s3]{Zhi-Xue Zhao}
\ead{zxzhao@amss.ac.cn}

\affiliation[s2]{
organization={School of Mathematical Sciences, Tianjin Normal University},
city={Tianjin},
postcode={300387}, 
country={China}
}

\affiliation[s3]{
organization={Institute of Mathematics and Interdisciplinary Sciences, Tianjin Normal University,},
city={Tianjin},
postcode={300387}, 
country={China}
}

\begin{abstract}
This study solves the output regulation problem for a reaction-diffusion system confronting concurrent input delay and fully unidentified disturbances (encompassing both unknown frequencies and amplitudes) across all channels. 
The principal innovation emerges from a novel adaptive control architecture that synergizes the modal decomposition technique with a dual-observer mechanism, enabling real-time concurrent estimation of unmeasurable system states and disturbances through a state observer and an adaptive disturbance estimator.
Unlike existing approaches limited to either delay compensation or partial disturbance rejection, our methodology overcomes the technical barrier of coordinating these two requirements through a rigorously constructed tracking-error-based controller, achieving exponential convergence of system output to reference signals. 
Numerical simulations are presented to validate the effectiveness of the proposed output feedback control strategy.
\end{abstract}


\end{frontmatter}

\section{Introduction}
\label{sec1}
\subsection{Background}

Output regulation is relevant in practical control scenarios, including the robot manipulation, the guidance of a missile towards a moving target, and the landing of an aircraft on a carrier (see \cite{1}). 
Specifically, output regulation aims to design a controller for an uncertain plant, wherein the plant is internally stable and the output asymptotically tracks the reference signal in the presence of disturbances. 
Output regulation problems for lumped parameter systems, described by ordinary differential equations (ODEs), have resulted in constructive outcomes and can be addressed by a set of linear matrix equations called the regulator equations (see \cite{14, 2}). 
On this basis, significant progress has been made in the output regulation for distributed parameter systems (see \cite{4,6,5}), including different types of parabolic partial differential equations (PDEs) (see \cite{13}).
For instance, a backstepping-based solution was proposed for the output regulation of coupled parabolic partial integro-differential equations (PIDEs) in \cite{7}, and observer-based controllers were developed for the output regulation of heat equations in \cite{9, 8}. 
As for the high-dimensional situation, Zhao {\emph{et al.}} derived an analytic feedback control for the output regulation of multi-dimensional heat equations by leveraging abstract results and PDE design techniques in \cite{11}. 
The previously mentioned research pertains to situations where the frequencies of disturbances are known. 
When frequencies are unknown, Guo {\emph{et al.}} utilized the adaptive observer internal model approach to estimate the frequencies for the output regulation of heat equations (see \cite{12,25}).

In this paper, we consider the output regulation for a reaction-diffusion equation.
Based on the regulator equations, we transform this problem into a stabilization problem for a reaction-diffusion equation with input delay and unknown frequency.
Since time delay has a detrimental impact on the stability of system, it is necessary to take it into account in control designs (see \cite{31}). 
Regarding stabilization problems for the reaction-diffusion equation with input delay, there are several effective approaches. 
In \cite{16}, \cite{17}, the backstepping-based transformations were utilized to map the reaction-diffusion equations into target systems, by which feedback controllers were subsequently developed to exponentially stabilize systems. 
Additionally, based on the feedback equivalence, Zhu {\emph{et al.}} proposed an observer-based feedback controller to exponentially stabilize the reaction-diffusion equation with input delay (see \cite{18}). 
In \cite{15,20}, the modal decomposition technique was combined with a predictor to compensate for input delay in reaction-diffusion equations. 
Unlike the former, Katz {\emph{et al.}} also employed modal decomposition technique, but combined it with LMIs to address the negative effects of input delay on these reaction-diffusion equations (see \cite{21,29}). 

Although there have been several results on the output regulation problem for reaction-diffusion equations (see \cite{22,23,26,10,24,34,27}), to the best of the authors know, the output regulation control design for this kind of model, which suffers from both input delay and disturbances with unknown frequency, has not been touched in the literature so far. 
In fact, the coexistence of these two challenges in the control system leads to significantly complex dynamics, thereby presenting substantial difficulties in controller design.

\subsection{Novelties}
This paper presents a novel solution to the output regulation problem for reaction-diffusion systems with dual challenges: unknown-frequency disturbances and input delay. 
We develop an innovative tracking-error-based controller that guarantees exponential reference tracking using only one single measurement.
 
1) We systematically investigate the critical impact of input delay on output regulation for reaction-diffusion PDEs. 
The proposed output feedback controller guarantees exponential decay of tracking errors, overcoming fundamental limitations in existing works.
While prior studies \cite{26,10,34,27} have advanced output regulation under time-varying disturbances, none address input delay effects - a pivotal practical consideration in distributed parameter systems.
When considering input delay, Gu {\emph{et al.}} required full-state feedback controllers for the output regulation in \cite{22,23}, and Wang imposed restrictive reaction coefficient constraints ($a \leq 1$) in \cite{24}. 
In contrast to these studies, we address a previously overlooked aspect-time delay-and propose a tracking error based controller, which is more implementable without restrictions on the reaction coefficient, that is, $a\in \mathbb{R}$.

2) This work tackles a more challenging scenario where disturbances with unknown frequency and amplitudes affect the entire system domain as well as both boundaries.
Unlike \cite{22,23,26,10,24} where disturbances only partially affect certain channels, we consider every channel of this reaction-diffusion system is subject to disturbances.
All existing works \cite{22,23,26,10,24,34,27} fundamentally rely on known disturbance frequencies - a restrictive assumption we eliminate. 
We develop an adaptive observer based on the adaptive internal model approach to dynamically track disturbances in real-time, which overcomes limitations of conventional Luenberger observers.

\subsection{Paper outline}

This paper is organized as follows.  
Section \ref{sec2} is divided into four parts: problem formulation, feedforward regulator design, error-based observer design, and error-based control design.
We give proof of most Lemmas, Propositions, and Theorems in Section \ref{sec3}.
Section \ref{sec4} and Section \ref{sec5} provide numerical simulations and conclusions, respectively.

\section{Statement of main results}
\label{sec2}

\subsection{Problem formulation}
The system we consider in this paper is described as 
\begin{equation}\label{1}
\left\{\begin{array}{ll}
w_t(x, t)=w_{x x}(x, t)+a w(x,t)+d_1(x) p(t),
\; x \in(0,1),  \\
w_x(0, t)=d_2 p(t),\\
w_x(1, t)=U(t-\tau)+d_3 p(t), \\
w(x, 0)=w_0(x),  \\
y_p(t)=w(0, t),
\end{array}\right.	
\end{equation}
where $y_p(t)$ is the output to be regulated, $U(t)$ is the control that satisfies $U(t-\tau) \equiv 0$ for $t-\tau<0$, and $\tau>0$ is the constant delay.
Additionally, $a\in\mathbb{R}$ is the reaction coefficient, and $w_0(x)$ is the initial state. 
We define $d_1(\cdot) \in L^{\infty}(0,1; \mathbb{R}^{1 \times 2})$, $ d_2 \in \mathbb{R}^{1 \times 2}$, and $d_3 \in \mathbb{R}^{1 \times 2}$ as unknown coefficients of the in-domain and boundary disturbances, respectively. 	
The disturbance $p(t)\in \mathbb{R}^{2\times1}$ and the reference signal $y_{\text{ref}}(t)$ are both generated by the following exosystem:
\begin{equation}\label{2}
\left\{\begin{array}{l}
\dot{p}(t)=G p(t),  \\
p(0)=p_0, \\
y_{\rm{ref}}(t)=d_4 p(t),
\end{array}\right.	
\end{equation}
where the initial value $p_0\in \mathbb{R}^{2 \times 1}$, $G \in \mathbb{R}^{2 \times 2}$, and $d_4\in\mathbb{R}^{1\times 2}$ are all unknown.
These mean that the reference signal $y_{\rm{ref}}(t)$ is also unknown.
We consider \eqref{1} and \eqref{2} in the state space $\mathcal{H}_1=L^2(0,1) \times \mathbb{R}^2$.

The following assumption is made throughout this paper.

\begin{assumption}\label{as1}
The spectrum of $G$ is $\{ i\omega, -i\omega\}$, where $\omega$ is an unknown positive parameter.		
\end{assumption}
 
By Assumption \ref{as1}, we know that \eqref{1} suffers from harmonic disturbances with unknown frequency $\omega$.	
Although the output $y_p(t)$ and the reference signal $y_{\rm{ref}}(t)$ are both unknown, the tracking error
$y_e(t):=y_p(t)-y_{\rm{ref}}(t)$ is known. 
The objective of the output regulation problem for \eqref{1} and \eqref{2} is to propose a tracking error based feedback controller such that 
\begin{equation}
\lim_{ t\to\infty }|y_e(t)| = 0, 
\; {\rm{exponentially}}.
\end{equation}	

\subsection{Feedforward regulator design}

In this section, we transform the output regulation of \eqref{1} and \eqref{2} into a stabilization problem by utilizing the regulator equations.
For this stabilization problem, we propose a feedforward regulator by the modal decomposition technique, LMIs, and Lyapunov-Krasovskii functionals.
First, we propose the following transformation:
\begin{equation}\label{86}
\varepsilon(x, t)=w(x,t)-\Gamma(x) p(t),
\end{equation} 
where $\Gamma: [0,1] \to \mathbb{R}^{1 \times 2}$ is the solution of the regulator equations
\begin{equation}\label{7}
\left\{\begin{array}{l}
\Gamma^{\prime \prime}(x)=\Gamma(x) G-a\Gamma(x)-d_1(x), \\
\Gamma^{\prime}(0)=d_2, \\
\Gamma(0)=d_4, 
\end{array}\right.	
\end{equation}
and it is easy to find the following Lemma.
\begin{lemma}\label{le7}
The initial value problem \eqref{7} has a unique solution $\Gamma(\cdot) \!\in\! H^2\big(\!(0,1);\!\mathbb{R}^{1 \times 2}\big)$.   
\end{lemma}

By transformations \eqref{86} and \eqref{7}, \eqref{1} is transformed into
\begin{equation}\label{8}
\left\{\begin{array}{l}
\varepsilon_t(x, t)=\varepsilon_{x x}(x, t)+a\varepsilon(x,t), \; t>0, \\
\varepsilon_x(0, t)=0, \\
\varepsilon_x(1, t)=U(t-\tau)-\gamma_1 p(t-\tau), \\
y_e(t)=\varepsilon(0, t),
\end{array}\right.	
\end{equation}
where $\gamma_1=(\Gamma^{\prime}(1)-d_3)e^{G\tau}$.

\begin{remark}\label{re2}
By Assumption \ref{as1}, we can extend the range of harmonic disturbance $p(t)$ to the interval $[-\tau,0]$.
As a result, we convert the output regulation problem of the closed-loop system $w(\cdot,t)\to\Gamma(\cdot)p(t)$ into a stabilization problem of the system $\varepsilon(\cdot,t)\to0$ in $L^2(0,1)$.
\end{remark}

Consider the Sturm-Liouville boundary value problem
\begin{equation}\label{103}
\left\{
\begin{array}{ll}
\phi^{\prime \prime}(x)+\lambda \phi(x)=0, \; 0<x<1, \\
\phi^{\prime}(0)=0,\\ 
\phi^{\prime}(1)=0,
\end{array}
\right.    
\end{equation}
where the eigenvalues $\{\lambda_n\}_{n=0}^{\infty}$ and the corresponding eigenfunctions $\{\phi_n\}_{n=0}^{\infty}$ are 
\begin{equation}\label{88}
\lambda_n=(n\pi)^2,
\;
\{\phi_n\}_{n=0}^{\infty}
=\left\{
\begin{array}{ll}
1, & n=0, \\
\sqrt{2}\cos n\pi x, & n\geq1. 
\end{array}
\right.    
\end{equation}
Therein, $\{\phi_n\}_{n=0}^{\infty}$ form an orthonormal basis of $L^2(0, 1)$, and we will use them to decompose variables in $L^2(0, 1)$.

Based on the modal decomposition technique, we propose the following feedforward regulator to stabilize \eqref{8}:
\begin{equation}\label{11}
U(t)=K\bar{\varepsilon}(t)+\gamma_{1} p(t),
\end{equation}
and
\begin{align} \label{17}
\bar{\varepsilon}(t)
&=[\langle\varepsilon(\cdot,t),\phi_0(\cdot)\rangle_{L^2(0,1)},
...,
\langle\varepsilon(\cdot,t),\phi_N(\cdot)\rangle_{L^2(0,1)}]^{\top} \cr
&= [\varepsilon_0(t),...,\varepsilon_N(t)]^{\top}
\in \mathbb{R}^{(N+1)}
\end{align}    
where $N$ will be determined by \eqref{18} later. 

\begin{remark}\label{re3}
The feedforward regulator \eqref{11} is composed of two parts: the first is to stabilize \eqref{8} without disturbance, and the second is to cancel the effect of disturbance. 
The feedback gain 
$
K=[k_0,k_1,...,k_N]\in\mathbb{R}^{1\times (N+1)}
$  
is determined by LMIs, for details, see Theorem \ref{le2}.
\end{remark}

Under the controller \eqref{11}, the system \eqref{8} reads
\begin{equation}\label{12}
\left\{
\begin{array}{l}
\varepsilon_t(x, t)=\varepsilon_{x x}(x, t)+a\varepsilon(x,t),\\
\varepsilon_x(0, t)=0, \\
\varepsilon_x(1, t)=
\left\{
\begin{array}{ll}
-\gamma_1p(t-\tau), & \tau>t>0, \\
K\bar{\varepsilon}(t-\tau), &t\geq\tau. \\
\end{array}
\right.
\end{array}
\right.	
\end{equation}

\begin{proposition}\label{le1}
For any initial value $\varepsilon(\cdot,0)\in L^2(0,1)$, under the control law \eqref{11}, there exists a unique weak solution $\varepsilon\in $
$C([0,\infty);L^2(0,1))$ to \eqref{12}.    
\end{proposition}

According to \eqref{103}, \eqref{88}, and Proposition \ref{le1}, the solution to \eqref{12} can be presented  by $\{\phi_n\}_{n=0}^{\infty}$ as
\begin{equation}\label{90}
\varepsilon(\cdot,t)
= \sum_{n=0}^{\infty} \varepsilon_n(t) \phi_n (\cdot),
\end{equation}
where 
\begin{equation}\label{95}
\varepsilon_n(t)
=\langle\varepsilon(\cdot,t),\phi_n\rangle_{L^2(0,1)}.  
\end{equation}
When $t>\tau$, we differentiate \eqref{95} along the solution of \eqref{12} and obtain
\begin{align}\label{14}
\dot{\varepsilon}_n(t)
& =\int^1_0 \varepsilon_{xx}(x,t)\phi_n(x) dx+a\int^1_0 \varepsilon(x,t)\phi_n(x) dx\cr
& =(-\lambda_n+a) \varepsilon_n(t) +K\phi_n(1)\bar{\varepsilon}(t-\tau).
\end{align}
Set $\delta>0$ as the desired decay rate of \eqref{12}, according to \eqref{88}, there exists $N\in\mathbb{N}$ such that 
\begin{equation}\label{18}
-\lambda_n+a<-\delta, \quad \forall n>N.
\end{equation}
Utilizing \eqref{18}, we separate \eqref{14} into the following two parts:
\begin{equation}\label{15}
\left\{\begin{array}{ll}
\dot{\bar{\varepsilon}}(t)=A\bar{\varepsilon}(t)+BK\bar{\varepsilon}(t-\tau), & n\leq N,\\
\dot{\varepsilon}_n(t)=(-\lambda_n+a)\varepsilon_n(t)+(-1)^n\sqrt{2}K\bar{\varepsilon}(t-\tau), & n> N,
\end{array}\right.	
\end{equation}    
where $\bar{\varepsilon}(t)$ is defined in \eqref{17} and 
\begin{equation}\label{89}
\left\{
\begin{array}{l}
A=\text{diag}(-\lambda_0+a,\ldots,-\lambda_N+a)\in\mathbb{R}^{(N+1)\times (N+1)},\\
B=[1\;-\!\sqrt{2}\;\ldots\;(-1)^N\!\cdot\!\sqrt{2}]^{\top}\in\mathbb{R}^{N+1},\\
K=[k_0\;k_1\;\ldots\;k_N]\in\mathbb{R}^{1\times (N+1)},
\end{array}
\right.
\end{equation}

\begin{theorem}\label{le2}
Suppose that $N$ satisfies \eqref{18} for any decay rate $\delta>0$, and 
$K=\left(k_0, k_1, \ldots, k_N\right)$. 
Let there exist matrices $P_2, P_3\in \mathbb{R}^{(N+1) \times (N+1)}$ and positive-definite matrices $P_1, S, R \in \mathbb{R}^{(N+1) \times (N+1)}$ such that 
\begin{equation}\label{97}
\Phi<0,
\end{equation}
where $\Phi=\{\Phi_{i j}\}$ is a symmetric matrix composed of
\begin{equation}\label{87}
\begin{array}{l}
\Phi_{11}=A^{\top} P_2+P_2^{\top} A+2 \delta P_1+S-e^{-2 \delta \tau} R, \\
\Phi_{12}=P_1-P_2^{\top}+A^{\top} P_3,
\quad 
\Phi_{13}=e^{-2 \delta \tau}R+P_2^{\top} BK, \\
\Phi_{22}=-P_3-P_3^{\top}+\tau^2 R,
\quad 
\Phi_{23}=P_3^{\top}BK,\\
\Phi_{33}=-e^{-2 \delta \tau}(S+R),
\end{array}
\end{equation}  
with $A,B$ given in \eqref{89}. 
For any initial value $\varepsilon(\cdot,0)\in L^2(0,1)$,  
there exists a unique solution $\varepsilon \in C([0,\infty)) ; L^2(0,1))$ to \eqref{12} such that
\begin{equation}\label{71}
\|\varepsilon(\cdot,t)\|_{L^2(0,1)}\leq M_1e^{-\delta t}\|\varepsilon(\cdot,0)\|_{L^2(0,1)},\;t\geq0,
\end{equation} 
for some $M_1>0$.
\end{theorem}
	
\subsection{Error-based observer design}
In this section, we first decouple the cascaded system \eqref{2} and \eqref{8} by a transformation. 
Then, for the PDE-part, we propose a state observer to track $\varepsilon(x,t)$ based on the tracking error $y_e(t)$. 
Finally, for the ODE-part, we design an adaptive observer to track the signal $\gamma_1p(t)$. 

According to Assumption \ref{as1}, $\gamma_1p(t)$ can be written as
$$
\gamma_1p(t)=b\cos\omega t+c\sin\omega t,
$$
where $b,c,\omega$ are unknown parameters and satisfy $b^2+c^2>0$.
Introducing a new variable 
$$
\eta(t)=[b\cos\omega t+c\sin\omega t \quad c\cos\omega t-b\sin\omega t]^\top,
$$ 
we rewrite $\gamma_1 p(t)$ as
\begin{equation}\label{41}
\left\{\begin{array}{l}
\dot{\eta}(t)=G_c \eta(t), \\
\gamma_1 p(t)=\gamma_c \eta(t),
\end{array}\right.
\end{equation}
where 
\begin{equation}
\gamma_c=[1,0], \; \eta(0)=(b, c)^{\top}, \;
G_c=\left[\begin{array}{cc}
0 & \omega \\
-\omega & 0
\end{array}\right].
\end{equation}	
Replacing \eqref{8}$_3$ with \eqref{41}$_2$, we obtain 
\begin{equation}\label{34}
\left\{\begin{array}{l}
\varepsilon_t(x, t)=\varepsilon_{x x}(x, t)+a\varepsilon(x,t), \\
\varepsilon_x(0, t)=0, \\
\varepsilon_x(1, t)=U(t-\tau)-\gamma_c \eta(t-\tau),\\
\dot{\eta}(t)=G_c \eta(t).
\end{array}\right.	
\end{equation}
Next, we design an observer for \eqref{34} only based on the tracking error $y_e(t)$.
Since \eqref{34} is ODE-PDE coupled, it is difficult to develop an observer for it directly.
To this end, we raise a transformation
\begin{equation}\label{30}
\epsilon(x,t)=\varepsilon(x, t)-g(x) \eta(t),
\end{equation}
and $g(x):[0,1] \to \mathbb{R}^{1\times2}$ satisfies
\begin{equation}\label{31}
\left\{\begin{array}{l}
g^{\prime \prime}(x)=g(x) (G_c-aE)-\mathcal{L}(x)g(0), \\
g^{\prime}(0)=0,\\
g^{\prime}(1)=-\gamma_c e^{-G_c\tau},
\end{array}\right.	
\end{equation}
where 
\begin{equation}\label{91}
\mathcal{L}(x)=\sum_{n=0}^{N}\ell_n\phi_n(x),
\end{equation} 
and $\ell_n$ is to be determined in Proposition \ref{th2}.

\begin{lemma}\label{le4}
The boundary value problem \eqref{31} admits a unique solution $g\in H^2\big((0,1);$ $\mathbb{R}^{1\times2}\big)$.
\end{lemma}

By the transformation \eqref{30}, we turn \eqref{34} into 
\begin{equation}\label{60}
\left\{\begin{array}{l}
\epsilon_t(x, t)=\epsilon_{x x}(x, t)+a\epsilon(x,t)+\mathcal{L}(x)(\epsilon(0,t)-y_e(t)),\\
\epsilon_x(0, t)=0, \\
\epsilon_x(1, t)=U(t-\tau),\\
\dot{\eta}(t)=G_c\eta(t).
\end{array}\right.	
\end{equation}	
Notably, \eqref{34} has been decoupled into two subsystems: `$\epsilon$'-part and `$\eta$’-part.
Next, we design a state observer for `$\epsilon$'-part
\begin{equation}\label{61}
\left\{\begin{array}{l}
\hat{\epsilon}_t(x, t)=\hat{\epsilon}_{x x}(x, t)+a\hat{\epsilon}(x,t)+\mathcal{L}(x)(\hat{\epsilon}(0,t)-y_e(t)),\\
\hat{\epsilon}_x(0, t)=0, \\
\hat{\epsilon}_x(1, t)=U(t-\tau),
\end{array}\right.	
\end{equation}	
which is a copy for `$\epsilon$'-part of \eqref{60}, but its initial condition is known.
Set the tracking error as
\begin{equation}\label{32}
\tilde{\epsilon}(x,t)=\epsilon(x,t)-\hat{\epsilon}(x,t),
\end{equation}
which satisfies
\begin{equation}\label{36}
\left\{\begin{array}{l}
\tilde{\epsilon}_t(x, t)=\tilde{\epsilon}_{x x}(x, t)+a \tilde{\epsilon}(x,t) +\mathcal{L}(x)\tilde{\epsilon}(0,t),\\
\tilde{\epsilon}_x(0, t)=0, \\
\tilde{\epsilon}_x(1, t)=0.
\end{array}\right.	
\end{equation}

By the classic semigroup theories, we can show the following result.
\begin{lemma}\label{le6}
For any initial value $\tilde{\epsilon}(\cdot,0)\in L^2(0,1)$, there exists a unique weak solution $\tilde{\epsilon}\in C([0,\infty);L^2(0,1))$ to \eqref{36}.	    
\end{lemma}     

According to Lemma \ref{le6}, the weak solution $\tilde{\epsilon}(\cdot,t)$ to \eqref{36} can be presented as
\begin{equation}\label{37}
{\tilde{\epsilon}}(x,t)=\sum_{n=0}^{\infty}{\tilde{\epsilon}}_n(t)\phi_n(x),	
\end{equation}
where 
\begin{equation}\label{96}
{\tilde{\epsilon}}_n(t)
=\langle {\tilde{\epsilon}}(\cdot,t),\phi_n \rangle_{L^2(0,1)}.
\end{equation}
Differentiate \eqref{96} along the solution of \eqref{36} yields
\begin{equation}\label{38}
\dot{{\tilde{\epsilon}}}_n(t)
= (-\lambda_n+a){\tilde{\epsilon}}_n(t)+{\tilde{\epsilon}}(0,t)\langle \mathcal{L}(\cdot),\phi_n\rangle_{L^2(0,1)}.
\end{equation}
According to $N$ defined in \eqref{18}, the system \eqref{38} can be separated into
\begin{equation}\label{39}
\left\{\begin{array}{ll}
\dot{\bar{e}}(t)=A\bar{e}(t)+LC\bar{e}(t)+L\left(\sqrt{2}\sum_{n=N+1}^{\infty}{\tilde{\epsilon}}_n(t)\right), 
& n\leq N,\\
\dot{{\tilde{\epsilon}}}_n(t)=(-\lambda_n+a){\tilde{\epsilon}}_n(t), 
& n> N,
\end{array}\right.	
\end{equation}    
where
\begin{equation}\label{33}
\left\{\begin{array}{l}	
\bar{e}(t)
=[{\tilde{\epsilon}}_0(t)\; \ldots \;{\tilde{\epsilon}}_N(t)]^{\top},
\;
L=[\ell_0\;\ldots\;\ell_N]^{\top},
\\
C=[1\;\sqrt{2}\;\sqrt{2}\;\ldots\;\sqrt{2}]\in\mathbb{R}^{1\times (N+1)},
\end{array}\right.
\end{equation}
and $A$ is defined in \eqref{89}.

\begin{proposition}\label{th2}
Suppose that $N$ satisfies \eqref{18} for the desired decay rate $\delta>0$, and $L=\left(\ell_0,\ell_1, \ldots, \ell_N\right)^{\top} \in \mathbb{R}^{N+1}$ can be chosen such that there exists a positive-definite matrix $Q\in \mathbb{R}^{(N+1)\times (N+1)}$ satisfying
\begin{equation}\label{50}
Q(A+LC)+(A+LC)^{\top}Q<-2\delta Q.
\end{equation}	
Then, for any initial value $\tilde{\epsilon}(\cdot, 0)\in L^2(0,1)$, there exists a unique solution 
$ \tilde{\epsilon} \in C([0,\infty]; L^2(0,1))$
to \eqref{36} such that
\begin{equation}
\|{\tilde{\epsilon}}(\cdot, t)\|_{L^2(0,1)} \leq M_2e^{-\delta t}\|{\tilde{\epsilon}}(\cdot, 0)\|_{L^2(0,1)}, \; t \geq 0,    
\end{equation}
for some $M_2>0$.
\end{proposition}

Next, we propose an adaptive observer to track $\gamma_c\eta(t)$. 
To this end, we set a new known variable 
\begin{equation}\label{92}
y_d(t)
:=y_e(t)-\hat{\epsilon}(0,t)
=g(0)\eta(t)+{\tilde{\epsilon}}(0,t),   
\end{equation}
which is based on the observer \eqref{61}. 
By \eqref{92}, the system \eqref{41} can be rewritten as
\begin{equation}\label{46}
\left\{
\begin{array}{l}
\dot{\eta}(t)=G_c\eta(t),\\
y_d(t)=g(0)\eta(t)+{\tilde{\epsilon}}(0,t).	
\end{array}
\right.	
\end{equation}

\begin{lemma}\label{le11}
The pair $(G_c,g(0))$ is observable for every $\omega\in(0,\infty)$.
\end{lemma}   
Lemma \ref{le11} guarantees that there exists an invertible transformation
\begin{equation}\label{57}
d(t)=T \eta(t)=(d_1(t), d_2(t))^{\top} \in \mathbb{R}^{2\times1},	
\end{equation} 
such that $d(t)$ satisfies the following canonical form
\begin{equation}\label{58}
\left\{\begin{array}{l}
\dot{d}(t)=S_c(\theta) d(t), \\
y_d(t)=\gamma_c d(t)+{\tilde{\epsilon}}(0, t),
\end{array}\right.	
\end{equation}
where	
\begin{equation}\label{83}
S_c(\theta)
=
\left[
\begin{array}{cc}
0 & 1 \\
-\theta & 0
\end{array}
\right], 
\theta=\omega^2,
\gamma_c=[1,0],     
\end{equation}
and $T$ is determined by
\begin{equation} \label{93}
 S_c(\theta)=T G_c T^{-1},
\gamma_c =g(0)T^{-1}.    
\end{equation}
In order to track $d(t)$ and $\theta$, we propose the following adaptive observer for \eqref{58}:		
\begin{equation}\label{54}
\left\{\begin{array}{l}
\dot{\xi}(t)=- \iota\xi(t)-y_d(t), \\	
\dot{\hat{\chi}}_1(t)=\hat{\varphi}(t)+\iota y_d(t)+\hat{\theta}(t) \xi(t)+k_0(y_d(t)-\hat{\chi}_1(t)), \\
\dot{\hat{\varphi}}(t)=-\iota \hat{\varphi}(t)-\iota^2 y_d(t), \\
\dot{\hat{\theta}}(t)=k_1 \xi(t)(y_d(t)-\hat{\chi}_1(t)),\\
\hat{d}_1(t)=\hat{\chi}_1(t), \\
\hat{d}_2(t)=\hat{\varphi}(t)+\xi(t) \hat{\theta}(t)+\iota \hat{\chi}_1(t),
\end{array}\right.   
\end{equation}
where $\iota>0, k_0>\frac{1}{4 \iota}, k_1>0$.

\begin{lemma}\label{le3}
For any initial value 
$(\xi(0), \hat{\chi}_1(0), \hat{\varphi}(0), \hat{\theta}(0))\in $ $ \mathbb{R}^4$, 
the observer \eqref{54} satisfies
$$
\lim _{t \rightarrow \infty}|\hat{\theta}(t)-\theta|=0 \; \rm{exponentially},
$$		
and 		
$$
\lim _{t \rightarrow \infty}\|\hat{d}(t)-d(t)\|_{\mathbb{R}^2}=0 \; \rm{exponentially}.	
$$
\end{lemma}

The proof of Lemma \ref{le3} is similar to claim 2 of \cite{33}, so we omit it. 
	
\subsection{Error-based control design}
In this section, we propose a tracking error based controller, which makes the output $y_p(t)$ track the reference signal $y_{\rm{ref}}(t)$ exponentially.

Set $f(x): [0,1] \to \mathbb{R}^{1 \times 2}$, which satisifes
\begin{equation}\label{48}
\left\{
\begin{array}{l}
f^{\prime \prime}(x)=f(x) S_c(\theta)-af(x)-\mathcal{L}(x)f(0), \\
f^{\prime}(0)=0, \\f(0)=\gamma_c.
\end{array}
\right.
\end{equation}	
Combining \eqref{31} and \eqref{93} with \eqref{48}, we obtain 
\begin{equation}\label{84}
f(x)=g(x)T^{-1}.    
\end{equation}
According to Lemma \ref{le4} and \eqref{84}, we have the following Lemma.
\begin{lemma}\label{le12}
The initial value problem \eqref{48} admits a unique solution 
$f(x) \in H^2\big( (0,1);$ $\mathbb{R}^{1 \times 2} \big)$, 
which is continuously differentiable concerning $x$ and the parameter $\theta$.
\end{lemma}

Together with \eqref{30}, \eqref{31}, \eqref{57}, \eqref{93} and \eqref{84}, we have 
\begin{equation}\label{49} 
\left\{\begin{array}{l}
\epsilon(x, t)
=\varepsilon(x, t)-f(x) d(t),\\
\gamma_c\eta(t)
=-g^{\prime}(1)e^{G_c\tau}\eta(t)
=-f^{\prime}(1)e^{S_c(\theta)\tau}d(t).
\end{array}	\right.
\end{equation}	
By \eqref{11}, \eqref{41}$_2$, \eqref{49} and observers \eqref{61} and \eqref{54}, we get the final feedback control law
\begin{align}\label{51}
U(t)
=& \sum_{n=0}^{N}k_n\langle\hat{\epsilon}(\cdot,t),\phi_n\rangle
+\sum_{n=0}^{N}k_n\langle f(\cdot,\hat{\theta}),\phi_n\rangle\hat{d}(t) 
-f^{\prime}(1,\hat{\theta})e^{S_c(\hat{\theta})\tau}\hat{d}(t).
\end{align}

\begin{lemma}\label{le5}
The tracking error based controller $U(t)$ defined in \eqref{51} converges exponentially to $K\bar{\varepsilon}(t)+\gamma_c\eta(t)$.
\end{lemma}
	
The final closed-loop system is described as follows:
\begin{equation}\label{94}
\left\{\begin{array}{ll}
w_t(x, t)=w_{x x}(x, t)+a w(x,t)+d_1(x) p(t), \\
w_x(0, t)=d_2 p(t),\\
w_x(1, t)=U(t-\tau)+d_3 p(t), \\
\dot{p}(t)=G p(t),  \\
y_e(t)
=w(0, t)-d_4 p(t),\\
U(t)
=\sum_{n=0}^{N}k_n\langle\hat{\epsilon}(\cdot,t),\phi_n\rangle
+\sum_{n=0}^{N}k_n\langle f(\cdot,\hat{\theta}),\phi_n\rangle\hat{d}(t) 
-f^{\prime}(1,\hat{\theta})e^{S_c(\hat{\theta})\tau}\hat{d}(t),\\
\hat{\epsilon}_t(x, t)=\hat{\epsilon}_{x x}(x, t)+a\hat{\epsilon}(x,t)+\mathcal{L}(x)(\hat{\epsilon}(0,t)-y_e(t)), \\
\hat{\epsilon}_x(0, t)=0, \;
\hat{\epsilon}_x(1, t)=U(t-\tau),\\
y_d(t)=y_e(t)-\hat{\epsilon}(0,t),\\
\dot{\xi}(t)=- \iota\xi(t)-y_d(t), \\	
\dot{\hat{\chi}}_1(t)=\hat{\varphi}(t)+\iota y_d(t)+\hat{\theta}(t) \xi(t)
+k_0(y_d(t)-\hat{\chi}_1(t)), \\
\dot{\hat{\varphi}}(t)=-\iota \hat{\varphi}(t)-\iota^2 y_d(t), \\
\dot{\hat{\theta}}(t)=k_1 \xi(t)(y_d(t)-\hat{\chi}_1(t)),\\
\hat{d}_1(t)=\hat{\chi}_1(t), \;
\hat{d}_2(t)=\hat{\varphi}(t)+\xi(t) \hat{\theta}(t)+\iota \hat{\chi}_1(t).
\end{array}\right.   	
\end{equation}

\begin{theorem}\label{th1}
Suppose that $N$ satisfies \eqref{18} for any decay rate $\delta>0$, and 
$ K=[k_0\;k_1\;\ldots\;k_N] $, 
$ L=[\ell_0\; \ell_1\; \ldots\; \ell_N] $ 
can be chosen such that there exists matrices $P_2, P_3\in \mathbb{R}^{(N+1) \times (N+1)}$, 
and
positive-definite matrices $P_1, S, R,Q \in$ $\mathbb{R}^{(N+1) \times (N+1)}$ 
satisfying \eqref{97} and \eqref{50}.
Then, for any initial value 
$( w(\cdot, 0), p(0), \hat{\epsilon}(\cdot, 0), \xi(0), \hat{\chi}_1(0), $
$ \hat{\varphi}(0), \hat{\theta}(0))
\in \mathcal{H}_1 \times L^2(0,1) \times \mathbb{R}^4$, 
the output tracking of \eqref{94} is guaranteed that
\begin{equation}\label{72}
\lim_{t\to\infty}|y_e(t)|=0,\;\rm{exponentially}.
\end{equation}
\end{theorem}

\section{Proof of main results}
\label{sec3}

\subsection{Proof of Proposition \ref{le1}}

According to Remark \ref{re2} and \cite[Theorem 1]{30}, we demonstrate that there exists a unique weak solution $\varepsilon\in C([0,\tau],L^2(0,1))$ to \eqref{12}.	
When $t\in[\tau, \infty)$, first, we set 
$
h(t):=K\bar{\varepsilon}(t)
$
and transfer $h(t)$ from the boundary of \eqref{12} into the PDE by invoking the following transformation 
\begin{equation}\label{65}
\tilde{\varepsilon}(x, t)=\varepsilon(x, t)-\frac{x^2}{2} h(t-\tau).
\end{equation}
By \eqref{65}, the system \eqref{12} is transformed into 
\begin{equation}\label{9}
\left\{\begin{array}{l}
\tilde{\varepsilon}_t(x, t)=\tilde{\varepsilon}_{x x}(x, t)+a \tilde{\varepsilon}(x, t)+f_1(x, t), \\
\tilde{\varepsilon}_x(0, t)=0, \quad \tilde{\varepsilon}_x(1, t)=0, \\
f_1(x, t)=\dfrac{ax^2+2}{2} h(t-\tau)-\dfrac{x^2}{2}\dot{h}(t-\tau) .
\end{array}\right.
\end{equation}	
The strategy for demonstrating the well-posedness of \eqref{12} on $[\tau,\infty)$
is to divide this interval by the delay $\tau$ and prove on each subinterval.
According to \cite[Theorem 1]{30} and \eqref{9}, when $t \in [\tau,2\tau]$, the proof for well-posedness of \eqref{12} reduces to proving 
$f_1 \in $ $ L^1([\tau,2\tau];L^2(0,1))$, 
which is sufficient to show $h(\cdot)\in C^1[0,\tau]$. 
To prove this, we utilize H\"{o}lder's inequality and have 
\begin{equation}\label{70}
|h(t_1)-h(t_2)|
\leq\sum_{n=0}^{N}|k_n|\|\varepsilon(\cdot,t_1)-\varepsilon(\cdot,t_2)\|_{L^2(0,1)},  
\end{equation}
and 
\begin{align}
 |\dot{h}(t_1)-\dot{h}(t_2)|
\leq & C_1\|\varepsilon(\cdot,t_1)-\varepsilon(\cdot,t_2)\|_{L^2(0,1)} \cr
& +\sqrt{2}\sum_{n=0}^{N}|k_n\gamma_1||p(t_1-\tau)-p(t_2-\tau)|,
\end{align}    
where $C_1=\sum_{n=0}^{N}|k_n|\left(\|\phi_n^{\prime\prime}\|_{L^2(0,1)}+|a|\right)<\infty$.
Since there exists a unique weak solution $\varepsilon\in C([0,\tau]; L^2(0,1))$ to \eqref{12} and $p(\cdot)\in C^1[-\tau,0]$ by Remark \ref{re2},
we obtain $h(\cdot)\in C^1[0,\tau]$, and further, 
$\varepsilon=\tilde{\varepsilon}+\frac{x^2}{2}h(\cdot-\tau) 
\in 
C([\tau, 2\tau]; L^2(0,1))$.
	
Then, we prove that \eqref{12} is well-posed on $[2\tau,3\tau]$. 
Since \eqref{12}$_3$ is different when $t>\tau$, this proof is somewhat different from the former subinterval. 
Specifically, we get the inequality
\begin{align}
|\dot{h}(t_1)-\dot{h}(t_2)|
\leq & C_1\|\varepsilon(\cdot,t_1)-\varepsilon(\cdot,t_2)\|_{L^2(0,1)} \cr
& +\sqrt{2}\sum_{n=0}^{N}|k_n||h(t_1-\tau)-h(t_2-\tau)|,
\end{align}    
which shows $h(\cdot)\in C^1[\tau,2\tau]$
and $\varepsilon \in C([2\tau, 3\tau]; L^2(0,1))$.
Using the proof process of $[2\tau,3\tau]$ on
$[m\tau,(m+1)\tau],\forall m=3,4,5,...$, 
with the initial condition $\varepsilon(\cdot,m\tau)\in L^2(0,1)$ 
and the continuous differentiability of $h(t-\tau)$ obtained on the previous interval, 
we finally get the existence of the unique weak solution $\varepsilon\in C([0,\infty); L^2(0,1))$.

\subsection{Proof of Theorem \ref{le2}}
By Parseval's identity and \eqref{90}, we obtain

\begin{equation}\label{79}
\|\varepsilon(\cdot,t)\|_{L^2(0,1)}^2
=\sum_{n=0}^N\varepsilon_n^2(t)+\sum_{n=N+1}^{\infty}\varepsilon_n^2(t).
\end{equation} 

Our first goal is to prove the first summand of \eqref{79} decays to zero exponentially. 
We propose the Lyapunov-Krasovskii functional 
$V_1=V_0+V_S+V_R$,
where
\begin{equation}\label{25}
\begin{aligned}
V_0(t) & ={\bar{\varepsilon}}^{\top}(t) P_1 \bar{\varepsilon}(t), \\
V_S(t) & =\int_{t-\tau}^t e^{-2 \delta(t-s)} {\bar{\varepsilon}}^{\top}(s) S \bar{\varepsilon}(s) d s, \\
V_R(t) & =\tau\int_{-\tau}^0 \int_{t+\theta}^t e^{-2 \delta(t-s)} \dot{\bar{\varepsilon}}^{\top}(s) R \dot{\bar{\varepsilon}}(s) ds d\theta,
\end{aligned}
\end{equation}
where the matrices $P_1,S,R$ are defined in Theorem \ref{le2}.
It is easy to find that there exists $C_2>0$ and $\kappa>0$ such that
\begin{equation}\label{81}
V_1(\tau)\leq C_2 e^{2\kappa\tau}|\bar{\varepsilon}(0)|^2\leq C_2 e^{2\kappa\tau}\|\varepsilon(\cdot,0)\|_{L^2(0,1)}^2.
\end{equation}
By  Jensen's inequality \cite[Proposition B.8]{31} and \eqref{15}$_1$, when $t \in [\tau,\infty)$, the derivative of three terms of $V_1(t)$ gives
\begin{equation}\label{19}
\begin{aligned}
\dot{V}_0+2 \delta V_0 
 = & 2 {\bar{\varepsilon}}^{\top} P_1 \dot{\bar{\varepsilon}}+2 \delta {\bar{\varepsilon}}^{\top} P_1 \bar{\varepsilon}, \\
\dot{V}_S+2 \delta V_S 
= & \bar{\varepsilon}^{\top} S \bar{\varepsilon}-e^{-2 \delta \tau}{\bar{\varepsilon}}^{\top}(t-\tau) S\bar{\varepsilon}(t-\tau), \\
\dot{V}_R+2 \delta V_R  
\leq &\tau^2 \dot{\bar{\varepsilon}}^\top R \dot{\bar{\varepsilon}} 
- e^{-2 \delta \tau}(\bar{\varepsilon}(t)-\bar{\varepsilon}(t-\tau))^{\top}R(\bar{\varepsilon}(t)-\bar{\varepsilon}(t-\tau)).
\end{aligned}
\end{equation}
It follows from \eqref{15}$_1$ that
\begin{equation}\label{21}
0=2[\bar{\varepsilon}^{\top} P_2^{\top}+\dot{\bar{\varepsilon}}^{\top} P_3^{\top}][-\dot{\bar{\varepsilon}}+A \bar{\varepsilon}+BK \bar{\varepsilon}(t-\tau)],
\end{equation}
where $P_2$ and $P_3$ are free-weighting matrices defined in Theorem \ref{le2}. 
Summing up the three terms of \eqref{19} and \eqref{21}, we obtain
\begin{equation}\label{26}
\dot{V}_1(t)+2 \delta V_1(t)\leq \psi^{\top}_1(t) \Phi \psi_1(t),
\end{equation}
where 
$
\psi_1=\operatorname{col}\{\bar{\varepsilon}(t), \dot{\bar{\varepsilon}}(t), \bar{\varepsilon}(t-\tau)\}
$ 
and $\Phi$ is defined in Theorem \ref{le2}.
Since $\Phi<0$, by \eqref{81} and \eqref{26}, there exists $C_3>0$ such that 
\begin{equation}\label{75}
\sum_{n=0}^N|\varepsilon_n(t)|^2
\leq\lambda_{\min}^{-1}(P_1)V_1(t)
\leq C_3e^{-2\delta t}\|\varepsilon(\cdot,0)\|_{L^2(0,1)}^2.
\end{equation}
Now, our problem reduces to proving that $\sum_{n=N+1}^\infty\varepsilon_n^2(t)$ converges to zero exponentially.
Solving \eqref{15}$_2$,	
by \eqref{18} and \eqref{75}, we find that there exists $C_4>0$ such that
\begin{equation}\label{80}
|\varepsilon_n(t)|
\leq e^{-\delta t}|\varepsilon_n(0)|+\frac{C_4\|\varepsilon(\cdot,0)\|_{L^2(0,1)}}{\lambda_n-a-\delta}e^{-\delta t}.
\end{equation}	
Since the series $\sum_{n=N+1}^{\infty}\limits\varepsilon_n^2(0)$, 
$\sum_{n=N+1}^{\infty}\limits\frac{C_4^2}{(\lambda_n-a-\delta)^2}$
both converge, by Minkowski inequality, \eqref{79}, \eqref{75} and \eqref{80}, we finally complete the proof. 

\subsection{Proof of Proposition \ref{th2}}
By Parseval's identity and \eqref{37}, when $t\in [\tau,\infty)$, ${\tilde{\epsilon}}(\cdot,t)$ can be represented as
\begin{equation}\label{44}	\|{\tilde{\epsilon}}\|_{L^2(0,1)}^2=\sum_{n=0}^{N}{\tilde{\epsilon}}_n^2(t)+\sum_{n=N+1}^{\infty}{\tilde{\epsilon}}_n^2(t).
\end{equation}
Our first goal is to demonstrate the second summand of \eqref{44} decays exponentially.
According to \eqref{18} and \eqref{39}$_2$, we have  
\begin{equation}\label{69}
\sum_{n=N+1}^{\infty}{\tilde{\epsilon}}^2_n(t)
\leq e^{-2\delta t}\left(\sum_{n=N+1}^{\infty}{\tilde{\epsilon}}^2_n(0)\right)
\leq e^{-2\delta t}\|{\tilde{\epsilon}}(\cdot,0)\|^2_{L^2(0,1)}.
\end{equation}
Regarding the first summand of \eqref{44}, we first need to demonstrate
$\sqrt{2}\sum_{n=N+1}^{\infty}{\tilde{\epsilon}}_n(t)$
in \eqref{39}$_1$ decays exponentially,
and we have the following equality:
\begin{equation}\label{67}
|\tilde{\epsilon}_n(t)|
= |{\tilde{\epsilon}}_n(0)|e^{(-\lambda_n+a)t}
\leq |{\tilde{\epsilon}}_n(0)|e^{(-\lambda_n+a)\tau}.
\end{equation} 
Combining H\"{o}lder's inequality and the Weierstrass M-test with \eqref{67}, we obtain the uniform convergence of 
$\sqrt{2}\sum_{n=N+1}^{\infty}{\tilde{\epsilon}}_n(t)$ and set
\begin{equation}\label{98}
\zeta(t)
:=\sqrt{2}\sum_{n=N+1}^{\infty}{\tilde{\epsilon}}_n(t)
:=\sqrt{2}e^{-\delta t}\zeta_0(t),
\end{equation}
where 
$
\zeta_0(t):=\sum_{n=N+1}^{\infty}{\tilde{\epsilon}}_n(0)e^{(-\lambda_n+a+\delta)t}
$. 
The task of proving $|\zeta(t)|$ decays exponentially reduces to proving that $\zeta_0(t)$ is uniformly bounded on $[\tau,\infty)$, and it suffices to prove $\zeta_0(t)$ decreases monotonically.
We get the following inequality:
\begin{equation}\label{64}
\sum_{n=N+1}^{m}|{\tilde{\epsilon}}_n(0)|e^{(-\lambda_n+a+\delta)t_2}\leq\sum_{n=N+1}^{m}|{\tilde{\epsilon}}_n(0)|e^{(-\lambda_n+a+\delta)t_1}.
\end{equation}
Due to the uniform convergence of $\zeta_0(t)$, we can take the limit respect to $m$ for both sides of \eqref{64}.
Consequently, we demonstrate that $\zeta_0(t)$ decreases monotonically, and further, $\zeta_0(t)$ is uniformly bounded.
Together with \eqref{98}, we have 
\begin{equation}\label{82}
|\zeta(t)|
\leq C_5 e^{-\delta t}\|\tilde{\epsilon}(\cdot,0)\|_{L^2(0,1)},
\end{equation}
for some $C_5>0$.
Regarding \eqref{39}$_1$, according to \eqref{89} and \eqref{33}, the pair $(A,C)$ is observable by the Hautus lemma. 
Thus, for the given $\delta$ and $N$ defined in \eqref{18}, $L$ can be chosen such that there exists a positive-definite matrix 
$
Q\in\mathbb{R}^{(N+1)\times (N+1)}
$ 
that satisfies the Lyapunov inequality \eqref{50}.   
On this basis, we propose a Lyapunov functional
\begin{equation}\label{99}
V_2(t)={\bar{e}}^{\top}(t) Q \bar{e}(t).    
\end{equation}
Differentiating \eqref{99} along the solution of \eqref{39}, by \eqref{50}, we have
\begin{equation}\label{104}
\dot{V}_2(t) +2\delta V_2(t)- \varpi_1 \zeta^2(t)
< \psi_2(t)^{\top} \Pi \psi_2(t),
\end{equation}
where 
$\psi_2(t) = {\rm{col}} \{\bar{e}(t) ,\zeta(t)\}$ 
and 
\begin{equation}
\Pi=
\left[
\begin{array}{cc}
-\varepsilon E  & QL \\
L^{\top} Q  & -\varpi_1
\end{array}
\right]
<0,
\end{equation}
for some  enough $\varepsilon>0$, some enough large $\varpi_1>0$, and $E\in \mathbb{R}^{(N+1) \times (N+1)}$.
By \eqref{104}, for an arbitrarily  $\varsigma>0$, we get
\begin{equation}\label{42}
\dot{V}_2(t) 
\leq -2(\delta+\varsigma) V_2(t)+\varpi_1\zeta^2(t) .
\end{equation}
Combining the comparison principle, \eqref{82}, and \eqref{42} with \eqref{44} and \eqref{69}, we finally complete the proof.
	
\subsection{Proof of Lemma \ref{le5}}
The error $\tilde{U}(t)$ between two control laws \eqref{11} and \eqref{51} is
\begin{align}\label{62}
&\tilde{U}(t) \cr
=& U(t)-(K\bar{\varepsilon}(t)+\gamma_c\eta(t))\cr
= & -\sum_{n=0}^{N}k_n\langle \tilde{\epsilon}(\cdot,t),\phi_n\rangle
+\sum_{n=0}^{N}k_n\langle f(\cdot,\hat{\theta})\hat{d}(t),\phi_n\rangle \cr
&-\sum_{n=0}^{N}k_n\langle f(\cdot,\theta)d(t),\phi_n\rangle 
+f^{\prime}(1,\hat{\theta})e^{S_c(\hat{\theta})\tau}\hat{d}(t)
-f^{\prime}(1,\theta)e^{S_c(\theta)\tau}d(t). 
\end{align}
According to Lemma \ref{le3}, $|\theta-\hat{\theta}(t)|$ converges exponentially to zero as $t \rightarrow+\infty$. 
Hence, we can set $\hat{\theta}(t), \theta \in[-C_6, C_6]$ for some $C_6>0$. 	
According to Lemma \ref{le12}, we have 
\begin{equation}\label{52}
\left\{\begin{array}{l}
\|f(x, \hat{\theta})-f(x, \theta)\|
\leq L_1|\hat{\theta}-\theta|
\leq L_1 C_6 e^{-\beta t}, \; x\in [0,1], \\
\|f^{\prime}(1,\hat{\theta})-f^{\prime}(1,\theta)\|
\leq L_1|\hat{\theta}-\theta|
\leq  L_1 C_6 e^{-\beta t},\\
\|e^{S_c(\hat{\theta})\tau}-e^{S_c(\theta)\tau}\|
\leq L_1|\hat{\theta}-\theta|
\leq  L_1 C_6 e^{-\beta t},
\end{array}\right.
\end{equation}
for some $L_1>0$. 
By \cite[Lemma 2.5]{12}, Lemma \ref{le3}, Proposition \ref{th2}, \eqref{62}, and \eqref{52},
we complete the proof.

\subsection{Proof of Theorem \ref{th1}}
Under the control law \eqref{51}, the system \eqref{8} becomes  	
\begin{equation}\label{53}
\left\{
\begin{array}{l}
\varepsilon_t(x, t)=\varepsilon_{x x}(x, t)+a\varepsilon(x,t) \\
\varepsilon_x(0, t)=0 \\
\varepsilon_x(1, t)=K\bar{\varepsilon}(t-\tau)+\tilde{U}(t-\tau).
\end{array}
\right.
\end{equation}
Differentiate \eqref{95} along the solution of \eqref{53} yields
\begin{equation}\label{55}
\dot{\varepsilon}_n(t)
= (-\lambda_n+a)\varepsilon_n(t)+(K\bar{\varepsilon}(t-\tau)+\tilde{U}(t-\tau))\phi_n(1).
\end{equation}
Separating \eqref{55} by $N$ defined in \eqref{18}, we obtain
\begin{equation}\label{56}
\left\{
\begin{array}{ll}
\dot{\bar{\varepsilon}}(t)=A\bar{\varepsilon}(t)+BK\bar{\varepsilon}(t-\tau)+B\tilde{U}(t-\tau), \; n\leq N,\\
\dot{\varepsilon}_n(t)=(-\lambda_n+a)\varepsilon_n(t)+(K\bar{\varepsilon}(t-\tau)+\tilde{U}(t-\tau))\phi_n(1),
n> N,
\end{array}
\right.	
\end{equation}     
where $A,B,K$ are defined in \eqref{89}.
Our first goal is to prove the first summand of \eqref{79} decays to zero exponentially. 
It follows from \eqref{56}$_1$ that
\begin{equation}\label{73}
0=
2[\bar{\varepsilon}^{\top} P_2^{\top}+\dot{\bar{\varepsilon}}^{\top} P_3^T]
[-\dot{\bar{\varepsilon}}
+A \bar{\varepsilon}
+BK \bar{\varepsilon}(t-\tau)
+B \tilde{U}(t-\tau)],
\end{equation}
where $P_2,P_3 \in \mathbb{R}^{(N+1) \times (N+1)}$ are two free-weighting matrices. 
Summing up the three terms of \eqref{19} and \eqref{73}, we have
\begin{equation}\label{102}
\dot{V}_\tau(t)+2\delta V_\tau(t)-\varpi_2 \tilde{U}^2(t-\tau)\leq \psi_3^\top \Lambda\psi_3,
\end{equation}
where
$\psi_3
=\operatorname{col}\left\{\bar{\varepsilon}(t), \dot{\bar{\varepsilon}}(t),\bar{\varepsilon}(t-\tau), \tilde{U}(t-\tau)\right\}$ 
and
\begin{equation}
\Lambda
=\left[
\begin{array}{ccc:c} 
&&& P_2^{\top} B \\
&\Phi& & P_3^{\top} B \\
&&& 0_{N \times 1} \\
\hdashline
B^{\top} P_2& B^{\top} P_3 & 0_{1 \times N} & -\varpi_2
\end{array}
\right]
<0,
\end{equation}
for some large enough $\varpi_2>0$ and $\Phi<0$ is defined in Theorem \ref{le2}.
By Lemma \ref{le5}, there exist $C_7>0$ and $\delta>\tilde{\omega}>0$ such that $\tilde{U}(t)\leq C_7 e^{-\tilde{\omega}t}$.
Utilizing the comparison principle to solve \eqref{102}, we can demonstrate the exponential decay of $\sum_{n=0}^{N}\varepsilon_n^2(t)$.

When $n>N$, solving \eqref{56}$_2$, by \eqref{18}, we have
\begin{equation} \label{101}
|\varepsilon_n(t)| 
\leq
|\varepsilon_n(0)|e^{-\delta t}
+\dfrac{C_8\|\varepsilon(\cdot,0)\|_{L^2(0,1)}}{\lambda_n-a-\tilde{\omega}}e^{-\tilde{\omega} t},
\end{equation}
for some $C_8>0$.
By H\"{o}lder's inequality, the exponential decay of $\sum_{n=0}^{N}\varepsilon_n^2(t)$, \eqref{101}, and
\begin{equation}
|y_e(t)| 
\leq \sqrt{2}  \sum_{n=0}^{N} |\varepsilon_n(t)|
+ \sqrt{2} \sum_{n=N+1}^{\infty} |\varepsilon_n(t)|, 
\end{equation}
we finally complete the proof.

\section{Numerical simulations}
\label{sec4}

\begin{figure}[htbp]
\centering
\subfigure[ Tracking performance ]
{
\includegraphics[width=2 in]{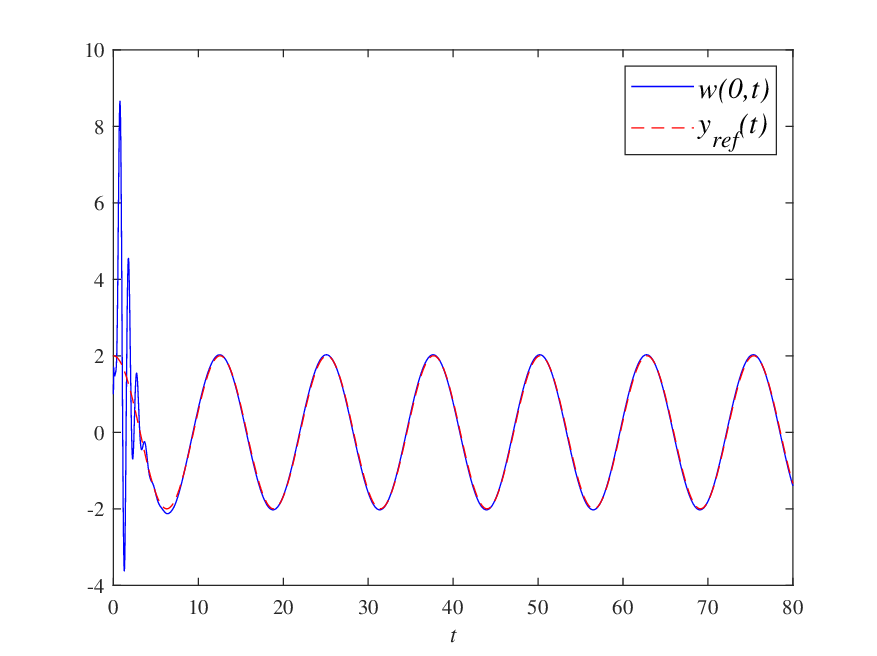}
\label{fig1}
} 
\quad 
\subfigure[ Evolution of $w(x,t)$ ]
{
\includegraphics[width=2 in]{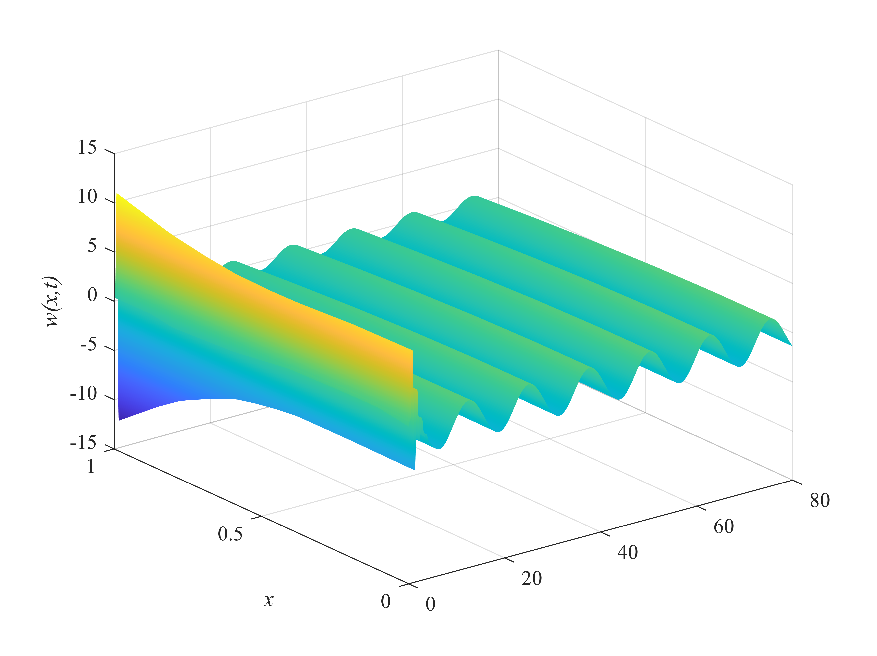}
\label{fig2}
}
\quad
\subfigure [ Estimate of $\theta$ ]
{
\includegraphics[width=4 in]{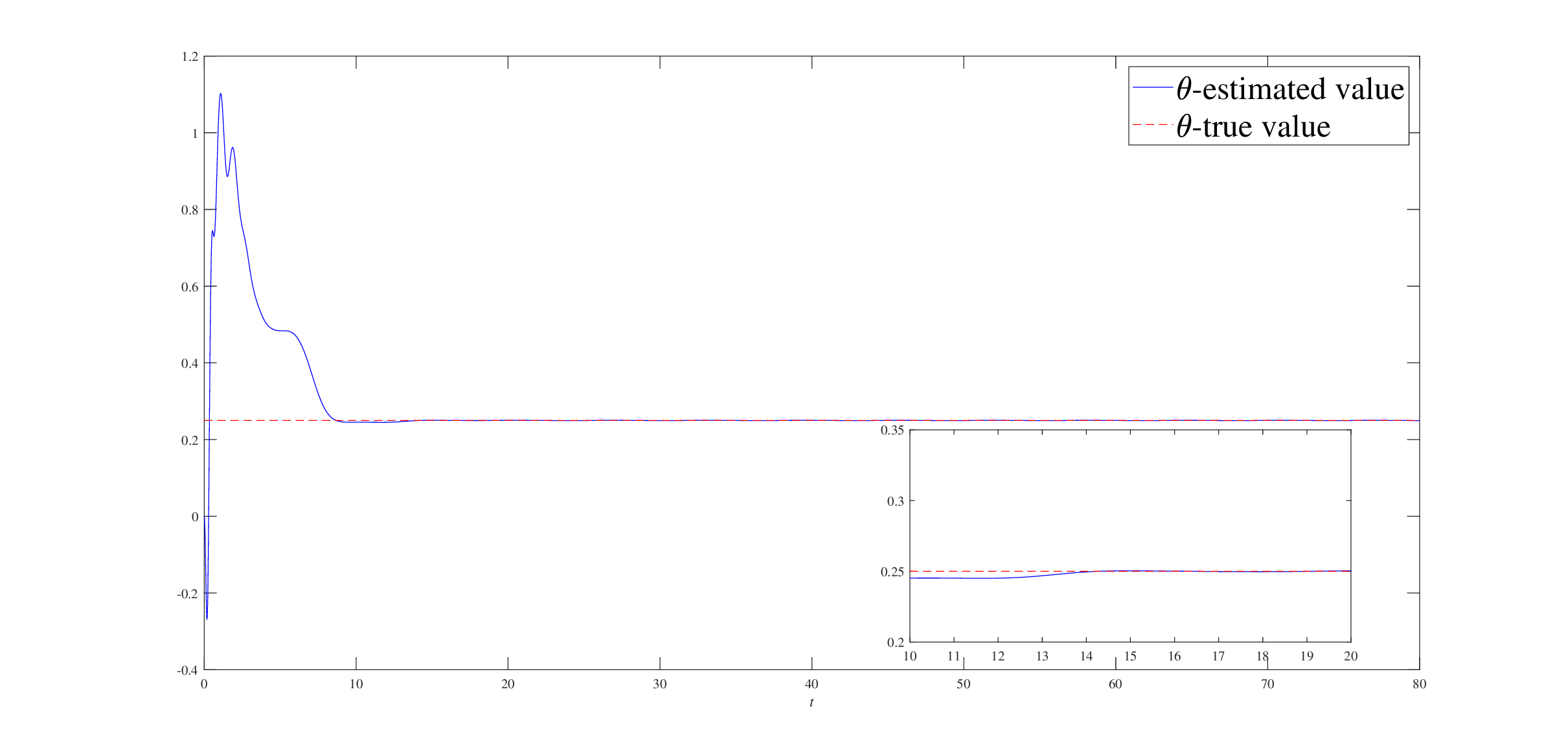}
\label{fig3}
} 
\caption{ Simulation results for $a=1.5$}
\label{fig7}
\end{figure}

\begin{figure}[htbp]
\centering
\subfigure[ Tracking performance ]
{
\includegraphics[width=2 in]{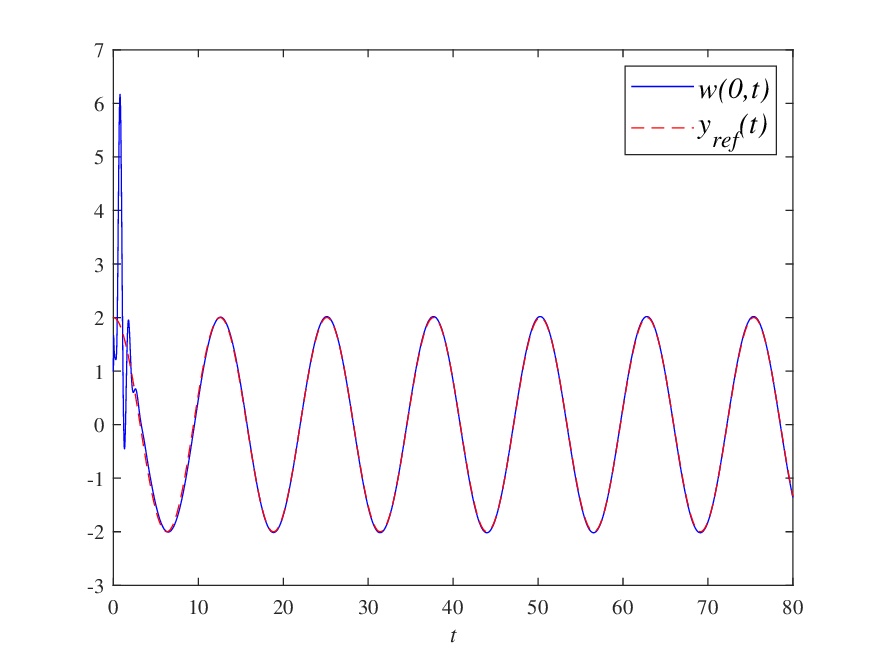}
\label{fig4}
} 
\quad 
\subfigure[Evolution of $w(x,t)$]
{
\includegraphics[width=2 in]{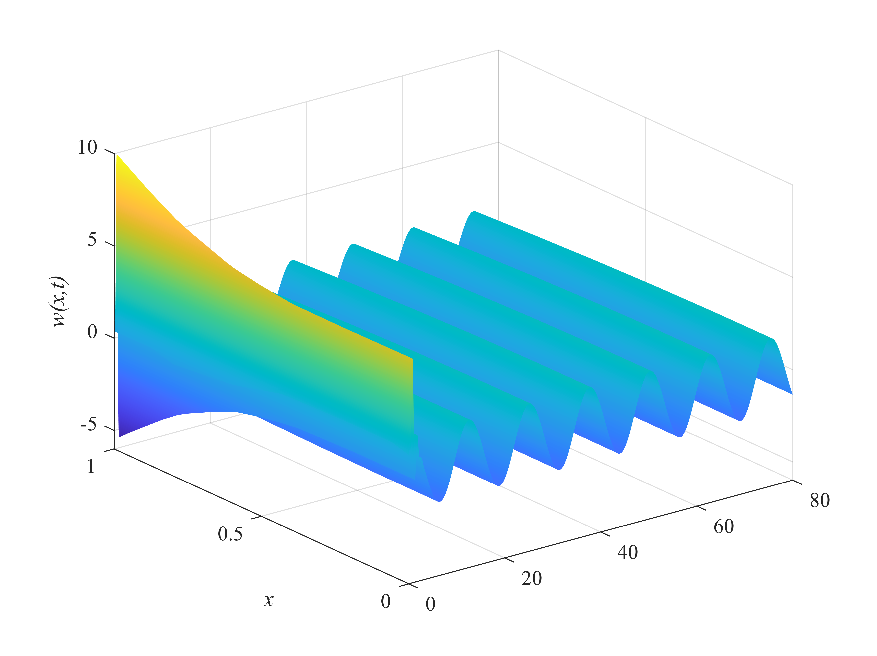}
\label{fig5}
}
\quad
\subfigure [Estimate of $\theta$]
{
\includegraphics[width=4 in]{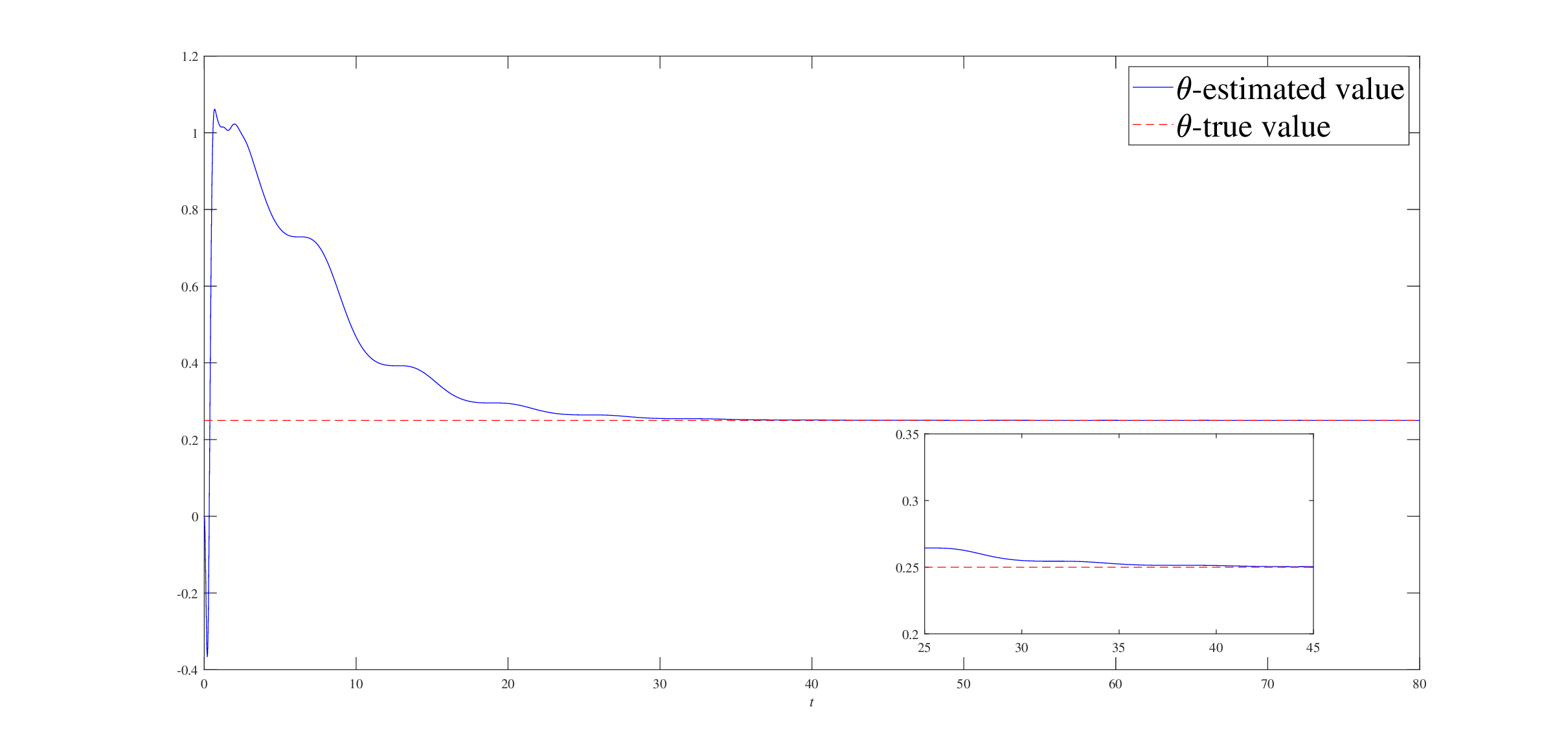}
\label{fig6}
} 
\caption{Simulation results for $a=0.5$}
\label{fig8}
\end{figure}

In this section, we present some numerical simulations to illustrate the effectiveness of the proposed controller \eqref{51}.
We consider the following system:
\begin{equation}\label{116}
\left
\{
\begin{array}{l}
w_{t}(x, t)=w_{x x}(x, t)+a w(x,t)+5x\sin(0.5t), \;x\in(0,1), \\
w_x(0, t)=\sin(0.5t), \\
w_x(1,t)=U(t-\tau),\\
y_p(t)=w(0,t),\\
y_e(t)=y_p(t)-y_{\rm{ref}}(t)=w(0,t)-2\cos(0.5t),
\end{array}
\right.
\end{equation}
The desired decay rate $\delta$ and parameters of the controller \eqref{51} are chosen as 
$
\delta=1,\; \tau=0.2,\;k_0=5, \; k_1=10, \; \iota=0.5.
$
According to \eqref{18}, we find the truncation order $N=1$. 
The initial values of the closed-loop system \eqref{94} are taken as
\begin{equation}\label{117}
\left\{
\begin{array}{l}
w(x,0)=\sin(2\pi x)+1,\;
\hat{\epsilon}(x, 0)=\sin(2\pi x)+1, \\
\xi(0)=\hat{\chi}_1(0)=\hat{\varphi}(0)=\hat{\theta}(0)=0.
\end{array}
\right.    
\end{equation}
The finite difference method is applied to numerically solve the system \eqref{94}. 
Fig. \ref{fig7} and Fig. \ref{fig8} describe the performance of \eqref{116} when the reaction coefficients are $a=1.5$ and $a=0.5$.

For Fig. \ref{fig7}, using the Matlab LMI toolbox, we find that the LMIs \eqref{87} and \eqref{50} in Theorem \ref{le2} and Proposition \ref{th2} have feasible solutions. 
By the resulting feasible solution, the control and observer gain matrices $K$ and $L$ in \eqref{89} and \eqref{33} are obtained as 
$
K=[-4.6591  \quad 0.2939], \; 
L= [-5.4933 \quad  7.7253]^\top.
$
Specifically, Fig. \ref{fig1} displays tracking performances of $y_p(t)$, which shows that $y_{p}(t)$ tracks $y_{\rm{ref}}(t)$ well.
Fig. \ref{fig3} indicates the performance of the frequency estimator of $\hat{\theta}(t)$ is satisfactory. 
Fig. \ref{fig2} demonstrate that the `$w$’-part is bounded.

It is worth noting that in Fig. \ref{fig8}, the reaction coefficient is $a=0.5$, which is smaller than 1, but the controller of \eqref{116} is still valid. 
The gains $K$ and $L$ are given as
$
K=[-3.6703 \quad 0.3204], \; 
L= [-4.8600 \quad 8.0813]^\top.
$
Fig. \ref{fig4} shows that $y_p(t)$ tracks $y_{\rm{ref}}(t)$ fast,
and Fig. \ref{fig6} demonstrates that the frequency estimator of $\hat{\theta}(t)$ performs well.
Fig. \ref{fig5} illustrates the boundedness of the `$w$’-part.

\section{Conclusion}
\label{sec5}

This paper presents a systematic solution to the output regulation problem for reaction-diffusion systems subject to input delay and full-channel disturbances with unknown frequency and amplitudes.
The considered system configuration is particularly general, requiring no restrictive assumptions on the reaction coefficient ($a\in \mathbb{R}$) or disturbance characteristics.
To address the input delay, we develop a compensation framework combining modal decomposition techniques with Lyapunov-Krasovskii functional analysis and LMIs.
For disturbance rejection, an adaptive observer incorporating the adaptive internal model principle is proposed to achieve real-time estimation of unknown frequency. 
The synthesized tracking error controller guarantees exponential convergence of the output to the reference signal.
Notably, the proposed control architecture exhibits extendability to broader classes of parabolic PDE systems beyond the considered reaction-diffusion equation.

The study also identifies two important open problems warranting future investigation. 
First, the extension to systems with multiple distinct input delays presents both theoretical and practical challenges, particularly regarding stability analysis under delayed actuation interactions. 
Second, the output regulation control design for cascaded PDE-ODE systems with time delays remains largely unexplored, representing a critical research direction. 
These unresolved problems highlight fundamental limitations in current distributed parameter system theory and suggest promising avenues for advancing the field.

\section*{CRediT authorship contribution statement}
\textbf{Shen Wang}: Writing – original draft, Investigation, Validation, Software, Methodology.
\textbf{Zhong-Jie Han}: Writing – review and editing, Supervision, Software, Methodology, Conceptualization.
\textbf{Kai Liu}: Writing – review and editing, Software, Methodology.
\textbf{Zhi-Xue Zhao}: Software, Methodology.
 
\section*{Declaration of competing interest}
The authors declare that they have no known competing financial interests or personal relationships that could have appeared to influence the work reported in this paper.

\section*{Acknowledgment}
This work was supported by the Natural Science Foundation of China grant NSFC-62473281, 12326327, 12326342. 

\appendix
\section{Proof of Lemma \ref{le4}}
We can find that $\varPsi_1=(1,i)^{\top}$ and $\varPsi_2=(1,-i)^{\top}$ are two eigenvectors of $G_c$ corresponding to the eigenvalues $\nu_1=i\omega$ and $\nu_2=-i\omega$, respectively.	
Right multiplying by $\varPsi_j,j=1,2$ to the ends of each equation in \eqref{31}, we work out the solution $g_j(x)=g(x)\varPsi_j \ne0,\;j=1,2$,
and further, obtain $g(x)$.

\section{Proof of Lemma \ref{le11} }
Since $(G_c,g(0))$ is observable if and only if $(J,\bar{g}(0))$ is observable, where $J={\rm{diag}}\{i\omega,-i\omega\},$ $\bar{g}(0)=(g_1(0)\;g_2(0))^\top$.
By Hautus lemma, we need $g_{1}(0)\ne0$ and $g_{2}(0)\ne0$, which is straight to see for $\omega\in(0,\infty)$ by proof of Lemma \ref{le4}. 

\bibliographystyle{elsarticle-num} 
\bibliography{reference.bib}

\end{document}